\documentclass[11pt, a4paper]{amsart}
\usepackage{amssymb}
\usepackage{amsthm}
\usepackage{amsmath, euscript}
\usepackage{color}
\usepackage{enumerate}
\usepackage[longnamesfirst]{natbib}
\usepackage[utf8]{inputenc}

\numberwithin{equation}{section}

\theoremstyle{plain}
\newtheorem{theorem}{Theorem}[section]

\newtheorem{cor}[theorem]{Corollary}
\newtheorem{prop}[theorem]{Proposition}

\theoremstyle{definition}
\newtheorem{definition}[theorem]{Definition}
\newtheorem{question}[theorem]{Question}

\theoremstyle{remark}

\numberwithin{equation}{section}
\newcommand{\until}{\upharpoonright}

\newcommand{\setN}{ \mathbb{N}}
\newcommand{\dom}{\mbox{dom}}

\newcommand{\la}{\langle}
\newcommand{\ra}{\rangle}
\newcommand{\wdt}[1]{\widetilde{#1}}
\newcommand{\Measures}{\EuScript{P}(2^\omega)}
\newcommand{\Csp}{{2^{\omega}}}
\newcommand{\Str}{{2^{<\omega}}}
\newcommand{\MLR}{\operatorname{MLR}}

\newcommand{\PA}{\mathrm{PA}}
\newcommand{\RR}{\mathcal{R}}
\newcommand{\lessT}{<_{\T}}
\newcommand{\leqT}{\leq_{\T}}
\newcommand{\geqT}{\geq_{\T}}
\newcommand{\eqT}{\equiv_{\T}}
\newcommand{\Tup}[1]{\ensuremath{\langle #1 \rangle}}

\DeclareMathOperator{\T}{T}

\newcommand{\rem}[1]{\relax}
\usepackage{ifthen, xspace}
    \newcommand{\bsc}{\usefont{T1}{cmr}{bx}{sc}}
    \newcommand{\pending}[1][]{{\noindent\highlight{\bsc\ifthenelse{\equal{#1}{}}{[to be written]}{[to be written: {\rm #1}]}}}\xspace}
    \newcommand{\prelims}[1][]{{\noindent\highlight{\bsc\ifthenelse{\equal{#1}{}}{[add to prelims]}{[add to prelims: {\rm #1}]}}}\xspace}
\usepackage{color}
    \newcommand{\highlight}[2][red]{{\color{#1}#2}}

\begin{document}

\title{Independence, Relative Randomness, and PA Degrees}
\date{\today}

\author[A.~R.~Day]{Adam R.~Day}
\address{Adam R.~Day\\
Department of Mathematic\\
University of California, Berkeley, CA \\
USA}
\email{adam.day@math.berkeley.edu}

\author[J.~Reimann]{Jan Reimann}
\address{Jan Reimann \\
Department of Mathematics\\
Pennsylvania State University\\
University Park, PA, USA}
\email{reimann@math.psu.edu}

\thanks{Reimann was partially supported by NSF grants DMS-0801270 and DMS-1201263.\\ \mbox{}}

\begin{abstract}
We study pairs of reals that are mutually Martin-L\"{o}f random with respect to a common, not necessarily computable probability measure. We show that a generalized version of van Lambalgen's Theorem holds for non-computable probability measures, too. We study, for a given real $A$, the \emph{independence spectrum} of $A$, the set of all $B$ so that there exists a probability measure $\mu$ so that $\mu\{A,B\} = 0$ and $(A,B)$ is $\mu\times\mu$-random. We prove that if $A$ is r.e., then no $\Delta^0_2$ set is in the independence spectrum of $A$. We obtain applications of this fact to PA degrees. In particular, we show that if $A$ is r.e.\ and $P$ is of PA degree so that $P \not\ge_{\T} A$, then $A \oplus P \geqT \emptyset'$.
\end{abstract}

\maketitle

\section{Independence and relative randomness}

The property of independence is central to probability theory. Given a probability space with measure $\mu$, we call two measurable sets  $\mathcal{A}$ and $\mathcal{B}$ independent if 
\begin{equation*}\mu \mathcal{A}= \frac{\mu (\mathcal{A} \cap \mathcal{B})}{\mu \mathcal{B}}.
\end{equation*}
The idea behind this definition  is that if event $\mathcal{B}$ occurs, it does not make event $\mathcal{A}$ any more or less likely. This paper considers a similar notion, that of relative randomness. The theory of algorithmic randomness provides a means of defining which elements of Cantor space ($2^\omega$) are random. We call $A \in 2^\omega$  Martin-L\"{o}f random if $A$ is not an element of any effective null set. We denote the class of all Martin-L\"{o}f random reals by $\MLR$.\footnote{For a comprehensive presentation of the theory of Martin-L\"{o}f randomness, see the monographs by \citet{downey-hirschfeldt:2010} and \citet{nies:2009a}. }

We say that $A$ is Martin-L\"{o}f random relative to $B$, or $A \in \MLR(B)$ if $A$ is not an element of any null set effective in $B$. Relative randomness is analogous  to independence because if $A \in \MLR(B)$, then not only is $A$ a random real but \textit{even given} the information in $B$, we cannot  capture $A$ in an effective null set. 
If we start with the assumption that $A$ and $B$ are both Martin-L\"{o}f random, then the following theorem of van Lambalgen establishes that relative randomness is symmetrical.

\begin{theorem}[\cite{VanLambalgen:1987}]
If $A,B \in \MLR$ then 
$A \in \MLR(B)$ if and only if $B \in \MLR(A)$ if and only if 
$A\oplus B \in \MLR$.  
\end{theorem}


We can extend the notion of relative randomness to any probability measure. We take $\Measures$ to be the set of all Borel probability measures on Cantor space. Endowed with the weak-$*$ topology, $\Measures$ becomes a compact metrizable space. The measures that are a finite, rational-valued, linear combination of Dirac measures form a countable dense subset, and one can choose a metric on $\Measures$ that is compatible with the weak-$*$ topology so that the distance between the those basic measures is a computable function, and with respect to which $\Measures$ is complete. In other words, $\Measures$ can be given the structure of an \emph{effective Polish space}. We can represent measures via Cauchy sequences of basic measures. This allows for coding measures as reals, and one can show that there exists a continuous mapping $\rho: \Csp \to \Measures$ so that for any $X \in \Csp$,
\[
	\rho^{-1}(\{\rho(X)\}) \text{ is a $\Pi^0_1(X)$ class.}
\]
For details of this argument, see \citet{day-miller:randomness-non-computable_2011}. If $\mu \in \Measures$, any real $R$ with $\rho(R) = \mu$ is called a \emph{representation} of $\mu$.

We want to define randomness relative to a parameter with respect to a probability measure $\mu$. Martin-Löf's framework easily generalizes to tests that have access to an oracle. However, our test should have access to \emph{two} sources: the parameter of relative randomness and the measure (in form of a representation).

\begin{definition}
Let $R_\mu$ be a representation of a measure $\mu$, and let $A \in \Csp$.

\begin{enumerate}[(a)]
 	\item A \emph{$(R_\mu,A)$-test} is given by a sequence $(\mathcal{V}_n \colon n \in \setN)$ of uniformly $\Sigma^0_1(R_\mu \oplus A)$-classes $\mathcal{V}_n \subseteq \Csp$ such that for all $n$, $\mu(\mathcal{V}_n) \leq 2^{-n}$.

 	\item A real $X \in \Csp$ \emph{passes} an $(R_\mu,A)$-test $(\mathcal{V}_n)$ if $X \not\in \bigcap_n \mathcal{V}_n$.

 	\item A real $X \in \Csp$ is $(R_\mu,A)$-\emph{random} if it passes all $(R_\mu,A)$-tests.
 
 \end{enumerate} 
\end{definition}

If, in the previous definition, $A = \emptyset$, we simply speak of an $R_\mu$\emph{-test} and of $X$ being \emph{$R_\mu$-random}. 

The previous definition defines randomness with respect to a specific representation. If $X$ is random for one representation, it is not necessarily random for other representations. 
On the other hand, we can ask whether a real exhibits randomness with respect to \emph{some} representation, so the following definition makes sense.
 
\begin{definition}
A real $X \in \Csp$ is \emph{$\mu$-random relative to $A \in \Csp$}, or simply \emph{$\mu$-$A$-random} if there exists a representation $R_\mu$ of $\mu$ so that $X$ is $(R_\mu,A)$-random. We denote by $\MLR_\mu(A)$ the set of all $\mu$-$A$-random reals. 
\end{definition}

For Lebesgue measure $\lambda$, we sometimes suppress the measure. Hence, in accordance with established notation, $\MLR(A)$ denotes the set of all  Martin-Löf random reals.

A most useful property of the theory of Martin-Löf randomness is the existence of \emph{universal tests}. Universal tests subsume all other tests. Furthermore, they can be defined uniformly with respect to any parameter. The construction can be extended to tests with respect to a measure $\mu$. More precisely, there exists a uniformly c.e.\ sequence $(U_n \colon n \in \setN)$ of sets $U_n \subseteq \Str$ such that, if we set for $R,A \in \Csp$
\[
	\mathcal{U}^{R,A}_n = \{ [\sigma] \colon \la \sigma, \tau_0, \tau_1 \ra \in U_n, \; \tau_0 \prec R, \tau_1 \prec A \},
\]
then $(\mathcal{U}^{R,A}_n)$ is an $(R,A)$-test and $X \in \Csp$ is $(R,A)$-random if and only if $X \not\in \bigcap_n \mathcal{U}^{R,A}_n$. We call $(U_n)$ a \emph{universal oracle test}.

Since for any $R \in \Csp$, $\rho^{-1}(\rho(R))$ is $\Pi^0_1(R)$, we can eliminate the representation of a measure in a test for randomness by defining, for any $A \in \Csp$,
\[
	\wdt{\mathcal{U}}^{R,A}_n = \bigcap_{S \in \rho^{-1}(\{\rho(R)\})} \mathcal{U}^{S,A}_n.
\]
The resulting class $\wdt{\mathcal{U}}^A_n$ is still $\Sigma^0_1(R)$, since $\rho^{-1}(\{\rho(R)\})$ is $\Pi^0_1(R)$ and hence compact. 

\begin{prop}
For any $R, A \in \Csp$ with $\rho(R) = \mu$, a real $X$ is $\mu$-$A$-random if and only if 
\[
 	X \not\in \bigcap_n \wdt{\mathcal{U}}^{R,A}_n. 
 \] 
\end{prop}

\begin{proof}
If $X$ is $\mu$-$A$-random, then it passes every $(R_\mu,A)$-test for some representation $R_\mu$ of $\mu$, in particular the instance $(\mathcal{U}^{R_\mu,A}_n)$ of the universal oracle test. Since $R_\mu \in \rho^{-1}(\{\rho(R)\})$, it follows that $X$ passes $\wdt{\mathcal{U}}^{R,A}_n$.

On the other hand, if for every representation $R_\mu$ of $\mu$, $X$ fails the test $(\mathcal{U}^{R,A}_n)$, then $X \in \bigcap_n \wdt{\mathcal{U}}^{R,A}_n.$
\end{proof}

The previous proposition shows that the test $\wdt{\mathcal{U}}^{R,A}_n$ is related to the concept of a \emph{uniform test}, originally introduced by \citet{Levin:1976vb}, and further developed by \citet{gacs:2005} and \citet{hoyrup-rojas:computabilityprobability_2009}. Hence we call it a \emph{uniform oracle test}.
Note that if $R,S$ are both representations of a measure $\mu$, then the uniform oracle tests $(\wdt{\mathcal{U}}^{R,A}_n)_n$ and $(\wdt{\mathcal{U}}^{S,A}_n)_n$ are identical.

\begin{definition}
Take $A, B \in 2^\omega$ and $\mu \in \Measures$. We say that $A$ and $B$ are \textit{relatively random with respect to $\mu$} if $A \in \MLR_\mu(B)$ and $B \in \MLR_\mu(A)$.
\end{definition}

Note that the representations of $\mu$ witnessing randomness for $A$ and $B$, respectively, do not have to be identical. 
If $A$ and $B$ are relatively random with respect to some measure $\mu$, then $\mu$ might offer some information about the relationship between $A$ and $B$. For example, we know that if $A$ and $B$ are relatively random with respect to Lebesegue measure, then any real they both compute must be K-trivial. If $A$ and $B$ are both atoms of $\mu$ then clearly $A$ and $B$ are relatively random with respect to $\mu$. Given this, perhaps the most obvious question to ask about relative randomness is the following. 

\begin{question}
For which $A, B \in 2^\omega$ does there exist a measure $\mu$ such that $A$ and $B$ are relatively random with respect to $\mu$ and neither $A$ nor $B$ is an atom of $\mu$?
\end{question}

This question is closely related to a theorem of \citet{Reimann:2008wy}. 
They proved that an element $X$ of Cantor space is non-recursive if and only if there exists a measure $\mu$ such that $X$ is $\mu$-random and $X$ is not an atom of $\mu$.

Van Lambalgen's theorem shows that $A$ and $B$ are relatively random if and only if $A\oplus B \in \MLR$. If we take $\lambda$ to be the uniform measure, then $A\oplus B \in \MLR$ if and only if the the pair $(A,B) \in 2^\omega \times 2^\omega$ is Martin-L\"{o}f random with respect to the product measure $\lambda \times \lambda$ i.e.\  $(A,B) \in \MLR_{\lambda \times \lambda}$. We begin our investigation into relative randomness  by showing that van Lambalgen's theorem holds for any Borel probability measure on Cantor space.

\begin{theorem} \label{thm:VL}
Let $\mu\in \Measures$ and let $A,B \in \Csp$ then
$(A,B) \in \MLR_{\mu \times \mu}$  if and only if $
A \in \MLR_\mu$  and $B \in \MLR_\mu(A)$. 
\end{theorem}

\begin{proof}
Let $R$ be any representation of $\mu$. 
First let us consider if $B \not \in \MLR_\mu(A)$. In this case we have that 
$B \in \bigcap_n \mathcal{U}_n^{R,A}$. We define an $(R,\emptyset)$-test for $\Csp \times \Csp$ by
$\mathcal{V}_n^R=\{[\tau]\times [\sigma] : \exists \eta \prec R \; (\Tup{\sigma,\eta,\tau} \in U_n)\}$. This ensures that $(A,B) \in \bigcap_n \mathcal{V}^R_n$. 
By applying Fubini's theorem we can establish that:
\begin{align*}
(\mu\times \mu) (\mathcal{V}^R_n) 
&=\int_{\Csp \times \Csp} \chi_{\mathcal{V}^R_n}(X,Y) d\mu \times d\mu\\
&=\int_\Csp \left ( \int_{\Csp } \chi_{\mathcal{U}_n^{R,X}}(Y) d\mu(Y) \right ) d\mu(X)\\
&\le \int_\Csp2^{-n}  d\mu(X) = 2^{-n}
\end{align*}
Hence $(A,B)$ is not $(R,\emptyset)$-random. As this is true for any representation $R$ of $\mu$ 
we have that $(A,B) \not \in \MLR_{\mu \times \mu}$. The same argument shows \textit{a fortiori} that if $A \not \in \MLR_\mu$ then $(A,B) \not \in \MLR_{\mu \times \mu}$.

To establish the other direction assume that 
$(A,B) \not \in \MLR_{\mu \times \mu}$. Again let $R$ be any
representation of $\mu$. 

Hence $(A,B) \in \bigcap_n \mathcal{V}^R_n$, where $(\mathcal{V}^R_n)$ is a universal $R$-test for $\Csp \times \Csp$. Let 
\[
	\mathcal{W}_n^{R,X}= \{Y : (X,Y) \in \mathcal{V}^R_n \}
\]
We have that
$\mathcal{W}_n^{R,X}$ is a $\Sigma^0_1(R\oplus X)$ class and this is uniform in $n$. However, given any $X$, we do not know whether or not  $\mu(\mathcal{W}_n^{R,X}) \le f(n)$ for some decreasing recursive function $f$  such that $\lim_nf(n) =0$. Hence we cannot necessarily turn this into a Martin-L\"{o}f test relative to $X$. In fact it is not even necessarily true that $\liminf_n \mu (\mathcal{W}_n^{R,X}) =0$.
We will show that the failure to turn this into a Martin-L\"{o}f test for some $X\in \Csp$ implies that $X \not \in \MLR_\mu$. 
This is a slight strengthening of the result that van Lambalgen obtained in his thesis. Van Lambalgen showed that if  $\liminf_n \mu (\mathcal{W}_n^{R,X}) \ne 0$ then $X \not \in \MLR_\mu$. 

However, we can generalize the proof of van Lambalgen's theorem given in \cite{nies:2009a}. We define another $R$-test by letting 
$\mathcal{T}_n^R = \{X \in \Csp: \mu(\mathcal{W}_{2n}^{R,X}) > 2^{-n}\}$.
To see that $\mathcal{T}_n^R \le 2^{-n}$, note that
\begin{align*}
(\mu\times \mu) \mathcal{V}^R_{2n}
&\ge\int_{\mathcal{T}_n^R \times \Csp } \chi_{\mathcal{V}^R_{2n}}(X,Y) d\mu \times d\mu\\
&=\int_{\mathcal{T}_n^R} \int_{\Csp } \chi_{\mathcal{W}_{2n}^{R,X}}(Y) d\mu(Y)  d\mu(X)\\
&\ge\int_{\mathcal{T}_n^R}2^{-n}  d\mu(X) = 2^{-n}\mu (\mathcal{T}_n^R)
\end{align*}
Now as $2^{-2n} \ge (\mu\times \mu) \mathcal{V}_{2n}^R$, we have that $\mu (\mathcal{T}_n^R) \le 2^{-n}$. Hence $\cap_n \mathcal{T}_n^R$ is an $R$-test.
Assume that $A \not \in \MLR_\mu(R)$. Then $A$ avoids all but finitely many of the sets $\mathcal{T}_n^R$. Hence for all but finitely many  $n$ we have that 
$\mu\, \mathcal{W}^{R,A}_{2n} \le 2^{-n}$ and so  by modifying finitely many $\mathcal{W}^{R,A}_{2n}$ we can obtain an $(R,A)$-test that covers $B$. Therefore $B$ is not $R$-random relative to $A$.

For all representations $R$ of $\mu$, we have shown that either $A$ is not $R$-random or $B$ is not $R$-random relative to $A$. However, to prove the theorem, it is essential that we get the \textit{same} outcome for all representations 
i.e.\ if $(A,B) \not \in \MLR_{\mu \time \mu}$ then either for all representations $R$ of $\mu$,  $A$ is not $R$-random or for all representations $R$ of $\mu$, $B$ is not $R$-random relative to $A$.

We can resolve this problem by taking our test $(\mathcal{V}^R_n)$ on the product space to be a uniform test. In this case we always obtain the same ``projection tests'' $(\mathcal{W}_n^{R,X})$ (independent of $R$) and hence the same outcome for any representation of $\mu$.
\end{proof}

\begin{cor}
If $A, B \in \Csp$ and $\mu \in \Measures$, then $A$ and $B$ are relatively random with respect to $\mu$ if and only if $(A,B) \in \MLR_{\mu \times \mu}$.
\end{cor}

\begin{cor}
\label{cor: not T comparable}
 If $A\ge_{\T} B$ and $(A,B) \in \MLR_{\mu \times \mu}$ then $B$ must be an atom of $\mu$.
\end{cor}

\begin{proof}
This holds because $B\in \MLR_\mu(A)$ if and only if $B$ is an atom of~$\mu$.
\end{proof}

We note that we cannot extend one direction of  van Lambalgen's theorem to  product measures of the form $\mu \times \nu$. In particular it is not true that if  $(A,B) \in \MLR_{\mu \times \nu}$ then $A \in \MLR_\mu$ and $B\in \MLR_\nu(A)$.  For example we can code $B$ into $\mu$ and obtain $A \in \MLR_\mu$, $B\in \MLR_\nu(A)$, $(A,B) \not \in \MLR_{\mu \times \nu}$.

\bigskip

Given any $X \in 2^\omega$, we will use $\RR(X)$ to denote the set of reals $Y$ such that $X$ and $Y$ are relatively random with respect to some measure $\mu$ and neither $X$ nor $Y$ are atoms of $\mu$ i.e.
\[
	\RR(A) = \{B\in 2^\omega: (\exists \mu \in \Measures)[(A,B) \in \MLR_{\mu \times \mu} \text{ and } \mu\{A,B\} =0 ]\}.
\]
We call $\RR(A)$ the \emph{independence spectrum of $A$}.

\medskip
The following proposition lists some basic properties of the independence spectrum.

\begin{prop}For all $A, B \in 2^\omega$ the following hold:
\begin{enumerate}
\item \label{p1} $A \in \RR(B)$ if and only if $B \in \RR(A)$. 

\item \label{p2} $B \in \RR(A) $ implies that $ A \mid_{\T} B$. 

\item \label{p3} If $A$ is non-recursive and $\nu$ is a computable, non-atomic measure (i.e.\ a measure with a computable representation and $\nu\{X\} = 0$ for all $X \in \Csp$), then $\RR(A)$ has $\nu$-measure $1$.

\item \label{p4} If $A \in \MLR$ then $\MLR(A) \subsetneq \RR(A)$. 
\end{enumerate}
\end{prop}
\begin{proof}
\eqref{p1} is by definition and  \eqref{p2} is by  Corollary~\ref{cor: not T comparable}.

\medskip
\eqref{p3} Suppose $A$ is non-recursive and $\nu$ is a computable measure with $\nu\{A\} = 0$. There is a measure $\mu$ such that $A$ is not an atom of $\mu$ and $A\in \MLR_\mu$, say via a representation $R_\mu$. Let $\kappa = (\mu + \nu)/2$. There exists a representation $R_\kappa \leqT R_\mu$, as $\nu$ is computable.	 
We claim that $A$ is $R_\kappa$-random. For if not, then $A$ fails some $R_\kappa$-test $(\mathcal{W}^{R_{\kappa}}_n)$. We have 
\[
	\mu \mathcal{W}^{R_{\kappa}}_n = 2\kappa \mathcal{W}^{R_{\kappa}}_n - \nu \mathcal{W}^{R_{\kappa}}_n \leq 2\kappa \mathcal{W}^{R_{\kappa}}_n \leq 2^{n-1}.
\]
Since $R_\kappa \leqT R_\mu$, $(\mathcal{W}^{R_{\kappa}}_{n+1})$ would define an $R_\mu$-test that covers $A$, contradicting the assumption that $A$ is $R_\mu$-random. Furthermore, by assumption on $\mu$ and $\nu$,  $\kappa\{A\}=0$. Hence 
\[
	(\MLR_\kappa(A) \setminus \{B : \kappa\{B\} \ne 0\}) \subseteq \RR(A),	
\] 
by van Lambalgen's Theorem. 

Now $\nu(\MLR_\kappa(A))=1$ because  the complement of $\MLR_\nu(A)$ is a $\kappa$ null set and hence a $\nu$ null set ($\nu$ is absolutely continuous with respect to $\kappa$ by definition). Moreover, the 
 set of atoms of $\kappa$ is countable and so has $\nu$-measure $0$ by the assumption that $\nu$ is non-atomic. This gives us that 
 \[
 	\nu(\MLR_\kappa(A) \setminus \{B : \kappa\{B\} \ne 0\}) = 1
 \]
 and thus $\nu \RR(A) = 1$. 

\medskip
\eqref{p4} Suppose $A$ is Martin-Löf random. By the definition of $\RR(A)$ and Theorem \ref{thm:VL} we have that $\MLR(A) \subseteq \RR(A)$. 

On the other hand, $A$ is not recursive and hence by \eqref{p3},
$\RR(A)$ has measure $1$ for any computable, non-atomic measure.
Let $\nu$ be a computable, non-atomic measure orthogonal to Lebesgue measure (e.g.\ the $(1/3,2/3)$-Bernoulli measure). Since $\nu \RR(A) = 1$, $\RR(A)$ has to contain a $\nu$-random element $X$. But $X$ cannot be relatively Martin-Löf random. Therefore, $\MLR(A) \subsetneq \RR(A)$.
%
%
\end{proof}

The proposition shows that, outside the upper and lower cone of a real $A$, the complement of $\RR(A)$ is rather small measure wise. On the other hand,
the above properties leave open the possibility that $\RR(A)$ is just the set of reals that are Turing incomparable with $A$. We will now establish that this is not necessarily the case.

\begin{prop}
Let $R$ be a representation of a measure $\mu$. If $A \in \Csp$ is such that
\begin{enumerate}
\item $A$ is r.e.\footnote{We mean here, of course, that $A$ is recursively enumerable viewed as a subset of $\setN$, by identifying a subset of $\setN$ with the real given by its characteristic sequence.},
\item $A$ is $R$-random, and 
\item $A$ is not an atom of $\mu$, 
\end{enumerate}
then $R\oplus A \ge_{\T} R'$.
\end{prop}

\begin{proof}
Given such an $R$ and $A$, let $A_s$  be a recursive approximation to $A$. We define the function $f\le_{\T} A \oplus R$ by:
\[f(x) = \min\{s :(\exists m \le s)(A_s\until m = A\until m \wedge
\mu_s[A\until m] <2^{-x})\}.\]
In this definition we take $\mu_s[\sigma]$ to be an $R$-recursive approximation to $\mu[\sigma]$ from above.
Note that $f$ is well defined because $A$ is not an atom of $\mu$. We claim that if $g$ is any partial function recursive in $R$, then for all but finitely many $x \in \dom(g)$, we have that $f(x) > g(x)$.
To establish this claim, let $g$ be an $R$-recursive partial function. We will build an $R$-test $\{U_n\}_{n\in \omega}$ by defining $U_n$ to be:
\[\{X \in \Csp : (\exists x>n) (\exists m)
 ( g(x)\downarrow \wedge \, \mu[A_{g(x)}\until m] <2^{-x} \wedge 
X \succ (A_{g(x)}\until m))\}.\]
Because any $x \in \dom(g)$ adds a single open set ($[A_{g(x)}\until m]$ for some $m$) of measure less than $2^{-x}$ to those $U_n$ with $n <x$, we have constructed a valid test.
Now if $g(x) \downarrow \ge f(x)$, then by definition of $f$, there is some $m \le f(x)$ such that
$\mu[A\until m] < 2^{-x}$ and $A\until m = A_{f(x)}\until m = A_{g(x)}\until m$. Thus for all $n < x$, $A \in U_n$.  Because $A \in \MLR_\mu(R)$ we have that $f(x) > g(x)$ for all but finitely many $x$ in $\dom(g)$. 

Let $g(x)$ be the $R$-recursive partial function with domain $R'$ such that $g(x)$ is the unique $s$ such that $x \in R'_{s+1}\setminus R'_s$. For almost all $x$, we have that $x\in R'$ if and only if $x\in R'_{f(x)}$ and so $R' \le_{\T}A \oplus R$.

\end{proof}

\begin{theorem}
\label{thm: measure and re}
Let $R$ be a representation of a measure $\mu$. If 
\begin{enumerate}
\item $A$ is r.e.,
\item $A$ is $\mu$-random, and 
\item $A$ is not an atom of $\mu$, 
\end{enumerate}
then $R\oplus A \ge_{\T} \emptyset'$.
\end{theorem}
\begin{proof}
Note the following characteristics of the previous proof. First the totality of $f$ does not depend on the fact that $A$ is $R$-random, it only depends on the fact that $A$ is not an atom of $\mu$. The construction is uniform so there is a single index $e$ such that  $\Phi_e(A\oplus \hat{R})$ is total if 
$\hat{R}$ is any representation of $\mu$. Additionally if $A$ is $\hat{R}$-random then for all but finitely many $x$, $\Phi_e(A\oplus \hat{R};x) \ge g(x)$ where $g$ is any $\hat{R}$-recursive partial computable function.

Let $R$ be any representation of $\mu$. The set 
$\{A \oplus \hat{R} : \hat{R}$ is a representation of $\mu\}$ is a
 $\Pi^0_1(A\oplus R)$ class and $\Phi_e$ is total on this class. From $A\oplus R$ we can compute a function $f$ that dominates $\Phi_e(A \oplus \hat{R})$ where $A$ is $\hat{R}$-random. As $f$ dominates any $\hat{R}$-recursive partial function we have that $A \oplus R \ge_{\T} \emptyset'$.
\end{proof}

\begin{cor}
If $A$ is r.e.\ and $B\le_{\T} \emptyset'$ then $B \not \in \RR(A)$. 
\end{cor}

The question remains, however, how big the independence of a real can be outside its upper and lower cones.

\begin{question}
	Is the set of all $X$ so that $X \mid_{\T} A$ and $X \not \in \RR(A)$ countable? 
\end{question}

\section{Recursively enumerable sets and $\PA$ degrees}

We will now give two (somehow unexpected) applications of Theorem~\ref{thm: measure and re} to the interaction between recursively enumerable sets and sets of $\PA$ degree. Recall that a set $A \subseteq \setN$ is of \emph{PA degree} if it is Turing equivalent to a set coding a complete extension of Peano arithmetic (PA). PA degrees have many interesting computability theoretic properties. For instance, a set is of PA degree if and only if it computes a path through every non-empty $\Pi^0_1$ class. However, a complete degree-theoretic characterization of the PA degrees is still not known. If $A \geqT \emptyset'$, then $A$ is of PA degree. On the other hand, Gödel's First Incompleteness Theorem implies that no r.e.\ set can be a complete extension of PA. \citet{jockusch-soare:degree_0} showed moreover that if a set is of incomplete r.e.\ degree, it cannot be of PA degree.  

It seems therefore worthwile to gain a complete understanding how r.e.\ sets and PA degrees are related. The crucial fact that links Theorem~\ref{thm: measure and re} to PA degrees is a result by \citet{day-miller:randomness-non-computable_2011}. They showed that every set of PA degree computes a representation of a \emph{neutral measure}.
Such a measure has the property that \emph{every} real is random with respect to it, i.e.\ $\Csp = \MLR_\mu$.  The existence of neutral measures was first established by \citet{Levin:1976vb}.

Our first result shows that below $\emptyset'$, r.e.\ sets and PA degrees behave quite complementary
\footnote{After the authors announced the result presented in Corollary \ref{cor:re-and-pa}, proofs not involving measure theoretic arguments have been found independently by A.~Ku{\v{c}}era and J.~Miller.}.

\begin{cor}[to Theorem~\ref{thm: measure and re}] \label{cor:re-and-pa}
If $A$ is an r.e.\ set and $P$ a set of  $\PA$ degree such that $P \not \ge_{\T} A$ then $P \oplus A \ge_{\T} \emptyset'$.
\end{cor}

\begin{proof}
  By the result of  \citet{day-miller:randomness-non-computable_2011} mentioned above, $P$ computes a representation $R_\mu$ of a neutral measure $\mu$ and $A \in \MLR_\mu$.  \citet{day-miller:randomness-non-computable_2011} also showed that a real $X$ is an atom of a neutral measure if every representation of the measure computes $X$. Now because $P \not \ge_{\T} A$, we have that $A$ is not an atom of $\mu$. Thus all hypotheses of Theorem \ref{thm: measure and re} are satisfied and we have $P \oplus A \ge_{\T} \emptyset'$.
\end{proof}

Corollary \ref{cor:re-and-pa} strengthens a result due to Ku{\v{c}}era and Slaman (unpublished). Recall that a function $f: \setN \to \setN$ is \emph{diagonally non-recursive} if $f(n) \neq \varphi_n(n)$ for all $n$, where $\varphi_n$ denotes, as usual, the $n$th partial recursive function.
Ku{\v{c}}era and Slaman constructed a a $\operatorname{low}_2$ r.e.\ set so that $A \oplus f \eqT \emptyset'$ for any diagonally non-recursive function $f \leqT \emptyset'$.
It is well-known that every PA degree computes a $\{0,1\}$-valued diagonally non-recursive function. Hence the set constructed by Ku{\v{c}}era and Slaman joins any PA degree below $\emptyset'$ to $\emptyset'$. Corollary \ref{cor:re-and-pa} yields that this is in fact true for \emph{any} r.e.\ set.

One can now ask  which kind of incomplete r.e.\ sets \emph{can} be bounded by PA degrees below $\emptyset'$.
This question was first raised by \citet{kucera:2004co}

\medskip
\begin{quote}
\it For which incomplete r.e.\ sets $A$ does there exist set $P$ of PA degree such that $A \lessT P \lessT \emptyset'$?
\end{quote}

\bigskip
We can use Corollary \ref{cor:re-and-pa} to completely answer this question.
We say a set $B$ is of \emph{PA degree relative} to a set $A$, written $B \gg A$ (see \citep{Simpson:1977ua}), if $B$ computes a path through every $\Pi^0_1(A)$ class. One well-known fact we will make use of is the following. If $P$ is of \emph{PA} degree, then there exists a set $Q$ of \emph{PA} degree such that $P \gg Q$. One way to prove this fact is to observe that the  $\Pi^0_1$ class
\[\{(A,B) \in 2^\omega \times 2^\omega \colon A  \in \mathrm{DNR}_2  \wedge
B\in \mbox{DNR}_2(A)\}\]
is non-empty, where $\mathrm{DNR}_2$ and $\mathrm{DNR}_2(A)$ are the classes of $\{0,1\}$-valued diagonally non-recursive functions and $\{0,1\}$-valued diagonally non-recursive functions relative to $A$, respectively.

\begin{theorem}
If $A$ is an r.e.\ set then the following are equivalent:
\begin{enumerate}
 \item $A$ is low.
 \label{it1}
\item There exists $P$, $P \gg A$ and $P$ is low.
\label{it2}
\item There exists $P$ of PA degree such that $ \emptyset' >_{\T} P>_{\T} A$.
\label{it3}
\end{enumerate}
\end{theorem}
\begin{proof}
(\ref{it1}) $\Rightarrow$ (\ref{it2}): There is a (non-empty) $\Pi^0_1(A)$ class of sets $B \gg A$. Relativize the low basis theorem to find $P \gg A$ and $P'\equiv_{\T} A'$. As $A$ is low so is $P$.

(\ref{it2}) $\Rightarrow$ (\ref{it3}): This is clear.

(\ref{it3}) $\Rightarrow$ (\ref{it1}): Take any  $Q$  of $\PA$ degree such that $P \gg Q$. 
Now $Q \ge_{\T} A$ because otherwise $Q \oplus A \ge \emptyset'$ but this is impossible because $P \ge_{\T} Q \oplus A$ and
$P \not \ge_{\T} \emptyset'$. Hence $P \gg A$. But now we have that $\emptyset'$ is r.e.\ in $A$ and also $\emptyset'$ computes a DNR function relative to $A$. Hence by relativizing Arslanov's completeness criterion we have that $A' \equiv_{\T} \emptyset'$.
\end{proof}

Observe that in the proof  of (\ref{it3}) $\Rightarrow$ (\ref{it1}),  showing $P \gg A$ only used the facts that $P \not \ge_T \emptyset'$ and $P\ge_T A$. Hence we get a final corollary.
\begin{cor}
If $P$ is a set of PA degree and $A$ is an r.e.\ set such that $P\ge_T A$ and $P \not \ge_T \emptyset'$, then 
$P \gg A$. 
\end{cor}

\section{Acknowledgments} 
\label{sec:acknowledgments}
The authors would like to thank Steve Simpson for stimulating and insightful discussions.


\bibliographystyle{abbrvnat}

\end{document}